\documentclass[11pt]{amsart}
\usepackage{amsfonts}
\usepackage{amsmath}
\usepackage{amssymb}
\usepackage{color}
\usepackage[all]{xy}
\usepackage{graphicx}
\usepackage{xspace}

 \font \eightrm=cmr8

 \newcommand{\nc}{\newcommand}

\newtheorem{thm}{Theorem}

\newtheorem{lem}[thm]{Lemma}
\newtheorem{prop}[thm]{Proposition}
\newtheorem{defn}{Definition}

\nc{\ignore}[1]{{}}
\nc{\mrm}[1]{{\rm #1}}
\nc{\dirlim}{\displaystyle{\lim_{\longrightarrow}}\,}
\nc{\invlim}{\displaystyle{\lim_{\longleftarrow}}\,}
\nc{\vep}{\varepsilon} \nc{\ep}{\epsilon}
\nc{\sigmat}{\widetilde\sigma}
\nc{\ostar}{\overline{*}}

\nc{\mchar}{\mrm{Char}}
\nc{\Hom}{\mrm{Hom}}
\nc{\id}{\operatorname{Id}}

\nc{\remark}{\noindent{\bf{Remark:}}}
\nc{\remarks}{\noindent{\bf{Remarks:}}}

 \nc{\delete}[1]{}
 \nc{\grad}[1]{^{({#1})}}
 \nc{\fil}[1]{_{#1}}

\nc{\BA}{{\Bbb A}} \nc{\CC}{{\Bbb C}} \nc{\DD}{{\Bbb D}}
\nc{\EE}{{\Bbb E}} \nc{\FF}{{\Bbb F}} \nc{\GG}{{\Bbb G}}
\nc{\HH}{{\Bbb H}} \nc{\LL}{{\Bbb L}} \nc{\NN}{{\Bbb N}}
\nc{\PP}{{\Bbb P}} \nc{\QQ}{{\Bbb Q}} \nc{\RR}{{\Bbb R}}
\nc{\TT}{{\Bbb T}} \nc{\VV}{{\Bbb V}} \nc{\ZZ}{{\Bbb Z}}

\nc{\Cal}[1]{{\mathcal {#1}}}
\nc{\mop}[1]{\mathop{\hbox {\rm #1} }}
\nc{\smop}[1]{\mathop{\hbox {\eightrm #1} }}
\nc{\mopl}[1]{\mathop{\hbox {\rm #1} }\limits}
\nc{\frakg}{{\frak g}}
\nc{\g}[1]{{\frak {#1}}}
\def \restr#1{\mathstrut_{\textstyle |}\raise-8pt\hbox{$\scriptstyle #1$}}
\def \srestr#1{\mathstrut_{\scriptstyle |}\hbox to
  -1.5pt{}\raise-4pt\hbox{$\scriptscriptstyle #1$}}
\nc{\wt}{\widetilde}
\nc{\wh}{\widehat}
\nc{\un}{\hbox{\bf 1}}
\nc{\redtext}[1]{\textcolor{red}{\tt #1}}
\nc{\bluetext}[1]{\textcolor{blue}{#1}}
\nc{\comment}[1]{[[{\tt {#1}}]] }
\nc{\R}{{\mathbb R}}

\nc\fleche[1]{\mathop{\hbox to #1 mm{\rightarrowfill}}\limits}
\def\semi{\mathrel{\times}\kern -.85pt\joinrel\mathrel{\raise 1.4pt\hbox{${\scriptscriptstyle |}$}}}

\begin{document}

\title[On enveloping algebras of preLie algebras]{\Large{Enveloping algebras of preLie algebras, Solomon idempotents and the Magnus formula}}

\author{Fr\'ed\'eric Chapoton}
\address{Institut Camille Jordan
Universit\'e Claude Bernard Lyon 1
B\^atiment Braconnier
21 Avenue Claude Bernard
F-69622 VILLEURBANNE Cedex
FRANCE}

\author{Fr\'ed\'eric Patras}
\address{Laboratoire J.-A.~Dieudonn\'e
         		UMR 6621, CNRS,
         		Parc Valrose,
         		06108 Nice Cedex 02, France}

\date{}

\begin{abstract}
We study the internal structure of enveloping algebras of preLie algebras.
We show in particular that the canonical projections arising from the Poincar\'e-Birkhoff-Witt theorem can be computed explicitely. They happen to be closely related to the Magnus formula for matrix differential equations. Indeed, we show that the Magnus formula provides a way to compute the canonical projection on the preLie algebra. Conversely, our results provide new insights on classical problems in the theory of differential equations and on recent advances in their combinatorial understanding.
\end{abstract}

\maketitle
\tableofcontents

\section{Introduction}
\label{sect:intro}
PreLie algebras play an increasingly important role in various fields of pure and applied mathematics, as well as in mathematical physics.
They appear naturally in the study of differential operators: the preLie product underlying the theory can be traced back to the work of Cayley
 in relation to tree expansions of compositions of differential operators. PreLie structures were defined formally much later in the work of Vinberg (they are sometimes refered to in the litterature as Vinberg algebras). They were rediscovered in deformation theory by Gerstenhaber \cite{gerst}, who introduced their current name.
They appear currently in differential geometry, control theory, operads, perturbative quantum field theory... see e.g. \cite{agra, ChaLiv01, EM} for some applications and further references on the subject.

In most applications, preLie algebras structures appear through their representations, 
and their action can therefore be canonically lifted to the enveloping algebra of the underlying Lie algebra.
The structure of enveloping algebras of preLie algebras, which are the subject of the present article, is richer than the structure of arbitrary enveloping algebras of Lie algebras. 
In the particular case of preLie algebras $L$, the usual Poincar\'e-Birkhoff-Witt (PBW) theorem can indeed be refined: 
there exists on $S(L)$, the symmetric algebra over $L$ another associative, simply defined, product, making $S(L)$ the enveloping algebra of $L$. Recall that, for Lie algebras, the PBW theorems asserts only that the canonical map from $S(L)$ to $U(L)$, the enveloping algebra constructed as the quotient of the tensor algebra by the ideal generated by Lie brackets, is a linear isomorphism; see e.g. \cite{r} for details on the subject.

In the context of differential operators and trees, the definition of the product on $S(L)$ first appeared in the work of Grossman-Larson \cite{gl}. 
The dual structure was discovered by Connes-Kreimer in the framework of perturbative quantum field theory \cite{Kreimer,CKII}.
The work of Chapoton-Livernet \cite{ChaLiv01} showed that free preLie algebras can be realized as Connes-Kreimer preLie algebras of trees, from which it follows that 
the enveloping algebra of a free preLie algebra $L$ can be realized as
a Grossman-Larson algebra, that is, as the symmetric algebra $S(L)$ provided with a new product associated to the process of insertion of trees. 
The general construction of the enveloping algebra of a preLie algebra $L$ as $S(L)$ is recent and due to Guin-Oudom, to which we refer for details \cite{GO}.

In the present paper, we study the internal structure of $S(L)$ in relation to the PBW theorem.
The classical combinatorial study of enveloping algebras, due originally to Solomon \cite{sol} and developped further in the work of Reutenauer and coauthors \cite{reutenauer1986,r}, shows that the canonical projections on the components of the direct sum decomposition of $U(L)$ induced by the PBW isomorphisms 
and the natural decomposition of the symmetric algebra into homogeneous components: $S(L)=\bigoplus\limits_{n\in \NN}S^n(L)$, 
where $L$ is the free Lie algebra over an arbitrary set of generators,
give rise to a family of orthogonal projections (i.e. idempotents) in the symmetric group group algebras $\QQ[S_n]$. 
These idempotents, the Solomon idempotents, play a key role in many fields of mathematics, partially due to their close connexion to enveloping algebra structures: 
free Lie algebras, control theory, group theory (both in the study of finite Coxeter groups and the general theory of transformation groups),
iterated integrals, Baker-Campbell-Hausdorff-type formulas (they were discovered independently of Solomon by Mielnik-Plebanski in this setting \cite{MP}), noncommutative symmetric functions, homological algebra, and so on... See e.g. \cite{r,g} for further details on the subject.
These idempotents have been given many names in the litterature: canonical idempotents (as canonical projections in the enveloping algebras), 
Barr idempotent (for the first idempotent of the series), eulerian idempotents (because, according to the seminal work of Solomon, 
their coefficients in the symmetric group algebras are related to the eulerian numbers). The one of ``Solomon idempotents'' fits the best with the practice for other Lie idempotents (e.g. the Dynkin, Garsia, Klyachko or Garsia-Reutenauer idempotents).

In the particular case of enveloping algebras of free preLie algebras, new structures and formulas arise in relation to the Solomon idempotents, due to the possibility of realizing these enveloping algebras as $S(L)$ equipped with a new product. This is the subject of the present article, which relates in particular the first Solomon idempotent to the Magnus series familiar in the theory of differential equations.

\section{PreLie and enveloping algebras}

In this article, all vector spaces are over an arbitrary ground field $k$ of characteristic zero.

\begin{defn}
 A preLie algebra is a vector space $L$ equipped with a bilinear map $\curvearrowleft$ such that, for all $x,y,z $ in $L$:
$$(x\curvearrowleft y)\curvearrowleft z-x\curvearrowleft (y\curvearrowleft z)=
(x\curvearrowleft z)\curvearrowleft y-x\curvearrowleft (z\curvearrowleft y).$$
\end{defn}
The vector space $L$ is then equipped with a Lie bracket $[x,y]:=x\curvearrowleft y-y\curvearrowleft x$.
We write $L_{Lie}$ for the associated Lie algebra, when we want to emphasize that we view $L$ as a Lie algebra. Most of the time this will be clear from the context, and we will then simply write $L$ although making use of the Lie algebra structure. For example, we will write $U(L)$ and not $U(L_{Lie})$ for the enveloping algebra of $L_{Lie}$. 
We will also denote by $\curvearrowleft$ the right action of the universal enveloping algebra of $L_{Lie}$ on $L$ that extends the pre-Lie product: $\forall a,b\in L, \ (b)a:=b\curvearrowleft a$. Notice that this action is well defined since the product $\curvearrowleft$ makes $L$ a module over the Lie algebra $L_{Lie}$:
$$\forall a,b,c\in L, \ ((c)b)a-((c)a)b=(c\curvearrowleft b)\curvearrowleft a-(c\curvearrowleft a)\curvearrowleft b$$
$$=c\curvearrowleft (b\curvearrowleft a- a\curvearrowleft b)=(c)[b,a].$$

Recall that one can consider $S(L)$, the symmetric algebra over $L$, equipped with a product $\ast$ induced by $\curvearrowleft$, as the enveloping algebra of $L_{Lie}$. Here, we realize $S(L)=\bigoplus\limits_{n\in\NN}S^n(L)$ as the subspace of symmetric tensors in $T(L):=\bigoplus\limits_{n\in\NN}L^{\otimes n}$, so that $S^n(L)$ identifies with the $S_n$-invariant tensors in $L^{\otimes n}$: $S^n(L)=(L^{\otimes n})^{S_n}$, where $S_n$ stands for the symmetric group of order $n$ -the permutation group of $[n]:=\{1,...,n\}$.
For $a\in L$, we write $a^n$ for its $n$-fold symmetric tensor power $a^{\otimes n}$ in $S(L)$. More generally, for $l_1,...,l_n\in L$, we write $l_1...l_n$ for the $n$-fold symmetric tensor product: $\frac{1}{n!}\sum_\sigma l_{\sigma(1)}\otimes ...\otimes l_{\sigma(n)}$, where $\sigma$ runs over $S_n$.

An enveloping algebra carries the structure of a cocommutative Hopf algebra for which the elements of the Lie algebra identify with the primitive elements. 
The corresponding coproduct on $S(L)$ is given by: for arbitrary $a_1,...,a_n\in L$,
$$\Delta(a_1...a_n)=\sum\limits_I a_I\otimes a_J$$
where, for a subset $I$ of $[n]$, $a_I:=\prod_{i\in I}a_i$, where $I$ runs over the (possibly empty) subsets of $[n]$ and $J:=[n]-I$.
The product $\ast$ is associative but not commutative and is defined as follows:
$$(a_1...a_l)\ast (b_1...b_m)=\sum\limits_f B_0(a_1\curvearrowleft B_1)...(a_l\curvearrowleft B_l),$$
where the sum is over all functions $f$ from $\{1,...,m\}$ to $\{0,...,l\}$ and $B_i:=f^{-1}(i)$. The increasing filtration of $S(L)$ by the degree is respected by the product $\ast$, but the direct sum decomposition into graded components is not. 

 A symbol $\,\widehat{}\,$ means that we consider the completion of a graded vector space: for example, $\widehat{L}$ stands for the product $\prod_n L_n$, where $L_n$ stands for the linear span of preLie products of length $n$ of elements of $L$. 
The series $\exp(a):=\sum\limits_n\frac{a^n}{n!}$ belongs for example to $\widehat{S(L)}$, and so on. A graded vector space $V=\bigoplus\limits_{n\in\NN}V_n$ is connected if $V_0=k$.

Since we are interested in universal properties of preLie algebras, we assume from now on that $L$ is a free preLie algebra. It follows that both $L$ and $L_{Lie}$ are naturally graded by the lengths of preLie products.
For clarity, we will distinguish notationally between:
 \begin{itemize}
  \item $U(L)$, the enveloping algebra of $L$ defined as the quotient of the tensor algebra $T(L)$ over $L$ by the ideal generated by the sums $[x,y]-x\otimes y+y\otimes x$, $x,y\in L$. The product on $U(L)$ is written simply $\cdot$ (so that e.g. for $x,y\in L, x\cdot y=x\otimes y$ in $U(L)$).
\item $S^\ast (L)$, the symmetric algebra over $L$ equipped with the structure of an enveloping algebra of $L_{Lie}$ by the product $\ast$. 
\item $S(L)$, the symmetric algebra over $L$ equipped with the structure of a polynomial algebra. Recall that the product of $l_1,...,l_n\in L$ in $S(L)$ is written $l_1l_2...l_n$.
 \end{itemize}
We write $S^{\ast,n}(L)$ and $S^n(L)$ for the symmetric tensors of length $n$ in $S^\ast(L)$ and $S(L)$; we write $S^{\geq n}(L)$ for sums of symmetric tensors of length at least $n$.
The coproduct on $S^\ast (L)$ and $S(L)$ for which $L$ is the primitive part is written $\Delta$ in both cases. It provides $S^\ast(L)$ with the Hopf algebra structure arising from its enveloping algebra structure, whereas $S(L)$ inherits a bicommutative Hopf algebra structure. Both $S^\ast(L)$ and $S(L)$ are free cocommutative coalgebras over $L$. We refer to \cite{mm, GO} for details.

The natural linear isomorphism (the PBW isomorphism) from $S(L)$ to $U(L)$ obtained by identifying $S(L)$ with the symmetric tensors in $T(X)$ is written $\pi_{nat}$:
$$\pi_{nat}(l_1...l_n)=\frac{1}{n!}\sum\limits_{\sigma\in S_n}l_{\sigma(1)}\cdot l_{\sigma(2)}\cdot ...\cdot l_{\sigma(n)}.$$ 
It is a coalgebra isomorphism since it maps $L$ to $L$ and since both $S(L)$ and $U(L)$ are cofree cocommutative coalgebras over $L$. Similarly, the identity map $\id$ from $S^\ast(L)$ to $S(L)$ is (obviously) a coalgebra isomorphism. 

The canonical isomorphism, written $is$, between $U(L)$ and $S^\ast(L)$ maps a product of elements of $L$ $x_1...x_n$ in $U(L)$ to their product $x_1\ast ...\ast x_n$ in $S^\ast (L)$. Take care that $is$ is not the inverse of the PBW map $\pi_{nat}$.

All these maps are summarized by the following diagram:
\begin{equation*}
  S^\ast(L) \stackrel{\id}{\longrightarrow} S(L) \stackrel{\pi_{nat}}{\longrightarrow} U(L) \stackrel{is}{\longrightarrow} S^\ast(L).
\end{equation*}

Recall at last from \cite{p0,p} the following two results:
\begin{lem}
 The Solomon first idempotent $sol_1$, that is the map from $U(L)$ (resp. $S^\ast(L)$) to $L$ orthogonally to the image of $S^{\geq 2}(L)$ by $\pi_{nat}$ (resp. $is\circ\pi_{nat}$) can be expressed as:
$$sol_1=\log^\star(\id),$$
where $\id$ is the identity map of $U(L)$, resp. $S^\ast(L)$.
More generally, the map to the image of $S^n(L)$ orthogonally to the image of the other graded components $S^i(L),\ i\not= n$, reads $$sol_n=\frac{sol_1^{\star n}}{n!}.$$
\end{lem}
The Lemma generalizes to arbitrary enveloping algebras of graded connected Lie algebras Reutenauer's computation of the Solomon idempotent in \cite{reutenauer1986}.
Here the logarithm is computed in the convolution algebras of linear endomorphisms of $U(L)$ and $S^\ast(L)$. That is, for $f,g\in End(U(L))$ (resp. $End(S^\ast(L))$):
$$f\star g:=m\circ (f\otimes g)\circ \Delta,$$
where $m$ stands for the product in $U(L)$, resp. $S^\ast(L)$.

\begin{prop}
 The inverse maps $inv$, resp $inv_S$ to the PBW isomorphisms $\pi_{nat}$ and $is\circ\pi_{nat}$ are given by:
$$inv(x):= \sum\limits_n\frac{sol_1^{\otimes n}}{n!}\circ\Delta^{[n]},$$
where $\Delta^{[n]}$ is the iterated coproduct map from $U(L)$ to $U(L)^{\otimes n}$. The same formula holds for $inv_S$.
\end{prop}

\section{The preLie PBW theorem}

Let us consider now the combinatorial PBW problem for preLie algebras, 
that is the particularization of Solomon's PBW problem for Lie algebras \cite{sol} to the preLie case:
compute explicit formulas for the canonical decomposition of the enveloping algebra $S^\ast(L)$ induced by the isomorphism $is\circ\pi_{nat}$ with $S(L)=\bigoplus\limits_nS^n(L)$.

We are interested in the general form of the decomposition, that is, the way an arbitrary element ${\underline l}=l_1...l_n$ of $S^\ast(L)$ decomposes
into ${\underline l}= sol_1({\underline l})+...+sol_n({\underline l})$. 
We can therefore assume that the $l_i$ are algebraically independent in $L$. Moreover, by the classical polarization argument, setting $a:=l_1+...+l_n$, we get that ${\underline l}$ is the $l_1,...,l_n$-multilinear component of $\frac{a^n}{n!}$; similarly, $sol_i({\underline l})$ identifies  to the $l_1,...,l_n$-multilinear component of $sol_i(\frac{a^n}{n!})$ and of $sol_i(\exp(a))$.
In the end, the computations can therefore be handled in the sub free preLie algebra of $L$ over $a$ and amount to computing the $sol_i(\exp(a))$. 

We get:
$$\exp(a)=\id (\exp(a))=(\exp^\star \circ\log^\star(\id)) (\exp(a)),$$
where we recall that $\log^\star(\id)=sol_1$.

Since $\Delta(a^n)=\sum\limits_{i\leq n}{n\choose i}a^i\otimes a^{n-i}$, 
${\underline a}:=\exp(a)$ is a group-like element in $S^\ast(L)$ (i.e. $\Delta ({\underline a})={\underline a}\otimes {\underline a})$.
However, the logarithm of a group-like element $x$ is a primitive element since 
$$\Delta (\log^\ast (x))=\log^\ast (x\otimes x)=
\log^\ast (x\otimes 1)\cdot (1\otimes x)=\log^\ast (x)\otimes 1+ 1\otimes \log^\ast (x)$$
and therefore $\log^\star(\id)({\underline a})=\log^\ast({\underline a})$ is a primitive element in $S^\ast(L)$ (here $\log^\star$, resp. $\log^\ast$ means that we compute the logarithm for the $\star$ product, resp. the $\ast$ product, and so on).
We get finally the first, polarized, form of the preLie PBW decomposition:
$${\underline a}=\sum\limits_n\frac{\pi^{[n]}\circ sol_1^{\otimes n}\circ\Delta^{[n]}}{n!}({\underline a})=\sum\limits_n\frac{(log^\ast({\underline a}))^{\ast n}}{n!},$$
where $\pi^{[n]}$ stands for the iterated $\ast$ product from $(S^\ast(L))^{\otimes n}$ to $S^\ast(L)$, resp. $\Delta^{[n]}$ for the iterated coproduct from $S^\ast(L)$ to    $(S^\ast(L))^{\otimes n}$.

Now, the $l_1,...,l_n$-multilinear part $sol_1(l_1...l_n)$ of $\log^\ast({\underline a})$ reads 
$$sol_1(l_1...l_n)=\sum\limits_{i=1}^n\frac{(-1)^{i-1}}{i}\sum_{I_1,...,I_i}l_{I_1}\ast...\ast l_{I_i}$$
where the $I_1,...,I_i$ run over all ordered partitions of $\{l_1,...,l_n \}$ ($I_1\coprod ...\coprod I_i=\{l_1,...,l_n\}$) where $\forall j\in [n],\ I_j\not=\emptyset$ and, for any subset $I$ of $[1,n]$,
$l_I:=\prod_{i\in I}l_i$. We get the first part of the following Proposition:

\begin{prop}
 The PBW decomposition of a generic element $l_1...l_n\in S^\ast (L)$ reads:
$$sol_1(l_1...l_n)=\sum\limits_{i=1}^n\frac{(-1)^{i-1}}{i}\sum_{I_1,...,I_i}l_{I_1}\ast...\ast l_{I_i}.$$
More generally, for higher components of the PBW decomposition, we have:
$$sol_i(l_1...l_n)=\sum\limits_{j=i}^n\frac{s(j,i)}{j!}\sum\limits_{I_1,...,I_j}l_{I_1}\ast ...\ast l_{I_j},$$
where the $s(j,i)$ are Stirling numbers of the first kind.
\end{prop}

The proof of the general case follows from a $\Psi$-ring argument as for the geometric computation of the $sol_i$ in \cite{pat91}.
Notice first that
$$\Psi^k:= \id^{\star k}=(\exp^\star (\log^\star (\id)))^{\star k}=\exp^\star (k\log^\star \id)=\sum_n k^n sol_n.$$
Besides, with the same notation as above for ordered partitions,  we get from the definition of the $\star$ product and the group-like behaviour of $\underline a$:
$$\Psi^{ k}(l_1...l_n)=\sum\limits_{j=1}^n\sum\limits_{I_1,...,I_j}{k\choose j}l_{I_1}\ast ...\ast l_{I_j}.$$
From the definition of the Stirling numbers $s(n,k)$ of the first kind, 
$$x(x-1)...(x-j+1)=\sum\limits_{1\leq i\leq j}s(j,i)x^j,$$
we get:
$$\Psi^{k}(l_1...l_n)=\sum\limits_{i=1}^n\left( \sum\limits_{j=i}^n\frac{s(j,i)}{j!} \sum\limits_{I_1,...,I_j}l_{I_1}\ast ...\ast l_{I_j}  \right) k^i$$
and finally, by identification of the coefficient of $k^i$:
$$sol_i(l_1...l_n)=\sum\limits_{j=i}^n\frac{s(j,i)}{j!}\sum\limits_{I_1,...,I_j}l_{I_1}\ast ...\ast l_{I_j},$$
which concludes the proof.

\section{The preLie PBW decomposition and the Magnus element}

\begin{defn}
  The Magnus element in the free preLie algebra $\widehat{L}$ over a single generator $a$ is
the (necessarily unique) solution $\Omega$ to the equation:
$$a\curvearrowleft \left(\frac{\Omega}{\exp(\Omega) -1}\right)=\Omega.$$
\end{defn}

The terminology is motivated by the so-called Magnus solution to an arbitrary matrix (or operator) differential equation 
$X'(t)=A(t)X(t),\ X(0)=1$:
$X(t)=\exp(\Omega(t))$, where
$$\Omega'(t)=\frac{ad_{\Omega(t)}}{\exp^{ad_{\Omega(t)}}-1}A(t)=A(t)+\sum\limits_{n>0}\frac{B_n}{n!}ad_{\Omega(t)}^n(A(t)),$$
where $ad$ stands for the adjoint representation and the $B_n$ for the Bernoulli numbers.
The link with preLie algebras follows from the observation that (under the hypothesis that the integrals and derivatives are well-defined), for arbitrary time-dependent operators, the preLie product
$$M(t)\curvearrowleft N(t):=\int_0^t [N(u),M'(u)]du$$
satisfies $(M(t)\curvearrowleft N(t))'=ad_{N(t)}M'(t)$. The Magnus formula rewrites therefore:
$$\Omega'(t)=\left(A(t)\curvearrowleft\left( \frac{\Omega}{\exp(\Omega)-1} \right)\right)'$$
where $\frac{\Omega}{\exp(\Omega)-1}$ is computed in the enveloping algebra of the preLie algebra of time-dependent operators. See e.g. the recent works by K. Ebrahimi-Fard and D. Manchon for further insights and an up-to-date point of view on the Magnus formula, in particular \cite{EM}.

\begin{thm}
 The Magnus element identifies with $sol_1(\exp(a))$:
$$\Omega = sol_1(\exp(a))=\log^\ast (\exp(a)).$$
\end{thm}

Indeed, we have:
$$sol_1(\exp(a))=\log^\ast({\underline a})=\sum_{n>0}(-1)^{n-1}\frac{({\underline a}-1)^{\ast n}}{n}.$$
Let us write $a=b+o(1)$ to mean that $a$ and $b$ in $S^\ast(L)$ are equal up to an element of $S^{\ast,\geq 2}(L)$.
We then have, by definition of the $\ast$ product: 
$$sol_1(\exp(a))= (a+o(1))\ast \left(\sum_n(-1)^{n-1}\frac{({\underline a}-1)^{\ast n-1}}{n}\right)$$
$$=a\ast  \left(\sum_n(-1)^{n-1}\frac{({\underline a}-1)^{\ast n-1}}{n}\right)+o(1).$$
Moreover, for an arbitrary element $b$ of $S^{\ast ,n}(L)$ the product $a\ast b$ reads $a\curvearrowleft b +o(1)$.
Finally, since $sol_1(\exp(a))\in L$, we get
$$sol_1(\exp(a))=a\curvearrowleft \left(\sum_n(-1)^{n-1}\frac{({\underline a}-1)^{\ast n-1}}{n}\right)=a\curvearrowleft\left(\frac{\log^\ast({\underline a})}{{\underline a}-1}\right),$$
from which the theorem follows.

\subsection*{Acknowledgements}
This work originated with discussions at the IESC conference Institute of Carg\`ese with K. Ebrahimi-Fard, F. Hivert, F. Menous, J.-Y. Thibon and the other participants to the CNRS PEPS program ``Mould calculus''. We thank them warmly for the stimulating exchanges, as well as the IESC 
for its hospitality.


\end{document}